\magnification=\magstep1
\centerline {\bf On Hecke $L$-functions attached to half-integral weight modular forms}
\bigskip\bigskip
\centerline {\it YoungJu Choie and Winfried Kohnen}

\bigskip\bigskip\bigskip
\noindent {\bf  }

We  investigate non-vanishing properties of $L(f,s)$ on the real line, when $f$ is a Hecke eigenform
  of half-integral weight $k+{1\over 2}$ on $\Gamma_0(4).$

\bigskip\bigskip\bigskip
\noindent {\bf 1. Introduction}
\bigskip
If $f$ is a cusp form of half-integral weight $k+{1\over 2}$ on $\Gamma_0(4)$ with Fourier coefficients $a(n)\, (n\geq 1)$, one can as usual attach to $f$ the Hecke $L$-series
$$L(f,s)=\sum_{n\geq 1}a(n)n^{-s}\quad (\sigma:=\Re(s)\gg1).\leqno (1)$$
It is known [9] that $L(f,s)$ has holomorphic continuation to ${\bf C}$ and satisfies a functional equation under $s\mapsto k+{1\over 2}-s$, relating $L(f,s)$ and $L(f|W_4,s)$ where $W_4$ is the Fricke involution.
\medskip
We note that even if $f$ is a Hecke eigenform, in general $L(f,s)$ has no Euler product. Though $L(f,s)$ has turned out to be useful in the investigation of sign changes of the coefficients $\Re(a(n))$ resp. $\Im(a(n))$ [6], in general otherwise its meaning remains a bit mysterious.
\medskip
In this paper we will investigate non-vanishing properties of $L(f,s)$ on the real line, when $f$ is a Hecke eigenform (Thm., sect. 3). The proofs which will be given in sect. 4 use the Mellin integral representation and rely on the existence of some special half-integral weight cusp forms with non-vanishing properties on the positive imaginary axis.
\medskip
As will be clear from the proof, similar (in fact, somewhat easier) arguments also work in the case of cuspidal Hecke eigenforms of integral weight on $\Gamma_1:=SL_2({\bf Z})$. We leave it to the reader to carry out the details and formulate the corresponding statements.
\medskip
We would like to recall that in the case of Hecke eigenforms on $\Gamma_1$ non-vanishing results for their Hecke $L$-functions at an {\it arbitrary} point $s_0$ in the critical strip (not on the critical line) have been proved in [4] (cf. also [7]), using holomorphic kernel functions. This method was carried over to the case of half-integral weight in [8], for arbitrary level. However, in this approach for given $s_0$ the weight in general has to be large depending on $s_0$.
\medskip
We also would like to point out that in [10] $L$-functions associated to modular forms of half-integral weight  were studied. In particular, certain half-integral weight cusp forms were investigated and the location of their zeros off the critical line were calculated.

\medskip\medskip
\leftline{{\bf 2010 Mathematics Subject Classification} Primary 11F37, Secondary 11F66 } 
\bigskip\bigskip

\medskip
Finally, addressing the reader interested in Siegel modular forms, we would like to mention the papers [3,5] where corresponding non-vanishing results for Koecher-Maass series are given.
\bigskip {\it Notations.} For $z\in {\bf C}^*$ we let $z^{{1\over 2}}:=e^{{1\over 2}\log z}$, where $\log z$ is the principal branch of the complex logarithm. For $z$ in the complex upper half-plane ${\cal H}$ we put $q=e^{2\pi iz}$.
\medskip
We denote by $\Gamma_0(4)$ the subgroup of $\Gamma_1$ consisting of matrices $\pmatrix{a&b\cr c&d\cr}$ with $4|c$.
\medskip
The letter $k$ always denotes an integer. We write $M_k(4)$ for the space of modular forms of weight $k$ for $\Gamma_0(4)$, with trivial character if $k$ is even and non-trivial character $({{-4}\over .})$ if $k$ is odd. We let $M_{k+{1\over 2}}(4)$ be the space of modular forms of weight $k+{1\over 2}$ for $\Gamma_0(4)$ and write
$S_{k+{1\over 2}}(4)$ for the subspace of cusp forms.
\bigskip
\noindent {\bf 2.  Basic facts on modular forms of half-integral weight of level 4}
\bigskip
For basic facts on modular forms of half-integral weight we refer to [9] and in the special case of level 4 also to [1,2].
\medskip
The group $\Gamma_0(4)$ has three cusps, represented by $0,{1\over 2}$ and $i\infty$. The cusp ${1\over 2}$ is $(k+{1\over 2})$-irregular, so a modular form of weight $k+{1\over 2}$ on $\Gamma_0(4)$ automatically vanishes at ${1\over 2}$.\medskip
As is well-known, one has
$$\dim S_{k+{1\over 2}}(4) =\sup\,\, \lbrace 0, [{k\over 2}]-1\rbrace.\leqno (2)$$
\medskip
Recall that on $M_{k+{1\over 2}}(4)$ (resp.  $S_{k+{1\over 2}}(4)$) the Fricke involution $W_4$ acts by
$$f(z)\mapsto (f|W_4)(z):=(-2iz)^{-k -1/2}f(-{1\over {4z}}) \quad (z\in {\cal H}).$$
The Hecke operators $T(p^2)$ ($p>2$ prime) commute with $W_4$, hence the spaces
$$S_{k+{1\over 2}}^{(\pm)}(4):=\lbrace f\in S_{k+{1\over 2}}(4)\,|\, f|W_4=(\pm)f\rbrace$$
have a basis of Hecke eigenforms of all the $T(p^2)$.
\medskip
For $f\in S_{k+{1\over 2}}(4)$ we put
$$L^*(f,s):=(2\pi)^{-s}\cdot 2^s\cdot \Gamma(s)\cdot L(f,s) \quad (\sigma\gg1)$$
where $L(f,s)$ is the $L$-series defined by (1). Then $L^*(f,s)$ has holomorphic continuation to ${\bf C}$ and satisfies the functional equation
$$L^*(f,k+{1\over 2}-s)=L^*(f|W_4,s).\leqno (3)$$
In particular, if $f|W_4=-f$, then $L^*(f,{k\over 2}+{1\over 4})=0$.
\medskip
In the rest of this section, we are concerned with certain special modular forms of level 4 which will play a role in  sect. 4.
\medskip
We let
$$\theta(z)=\sum_{n\in {\bf Z}}q^{n^2} \quad (z\in {\cal H})$$
be the basic theta function which is in $M_{{1\over 2}}(4)$. One has
$$\theta|W_4=\theta. \leqno (4)$$
\medskip
Let
$$P(z)=1-24\sum_{n\geq 1}\sigma_1(n)q^n  \quad (z\in {\cal H})$$
be the quasi-Eisenstein series of weight 2 on $\Gamma_1$ (where $\sigma_1(n)=\sum_{d|n}d$) and put
$$F_2(z):={1\over 24}\bigl(-P(z)+3P(2z)-2P(4z)\bigr).$$
Then
$$F_2(z)=\sum_{n\geq 1,n\equiv 1\pmod 2}\sigma_1(n)q^n$$
and $F_2\in M_2(4)$. The transformation formula
$$(2z)^{-2}F_2(-{1\over {4z}})=F_2(z)-{1\over {16}}\theta^4(z)\leqno (5)$$
holds.
\medskip
Further, put
$$\Delta_4(z):= F_2(z)\bigl(\theta^4(z)-16F_2(z)\bigr) \quad (z\in {\cal H}).$$
Then $\Delta_4\in M_4(4)$ and one has
$$(2z)^{-4}\Delta_4(-{1\over {4z}})=\Delta_4(z)\leqno (6)$$
as follows form (4) and (5).
\medskip
We also note that $\Delta_4$ vanishes at the cusps $0$ and $i\infty$, hence it follows that $\theta\Delta_4$ is a cusp form of weight ${9\over 2}$ on $\Gamma_0(4)$. Using (2) we find that the map
$$f\mapsto f\cdot \theta\Delta_4$$
gives an isomorphism of $M_k(4)$ onto $S_{k+{9\over 2}}(4)$.
\medskip
Finally, we define
$$D_2(z):=\theta^4(z)-32 F_2(z) \quad (z\in {\cal H}).$$
Then $D_2$ is in $M_2(4)$ and using (4) and (5) one checks that
$$(2z)^{-2}D_2(-{1\over {4z}})=D_2(z).\leqno (7)$$
\bigskip
\noindent{\bf 3. Statement of results}
\bigskip
It easily follows from the discussions in sect. 2 that $S_{k+{1\over 2}}^{(+)}(4)=\lbrace 0\rbrace$ for $k<4$ (in fact, $S_{k+{1\over 2}}(4)=\lbrace 0\rbrace$ in this case) and  $S_{k+{1\over 2}}^{(-)}(4)=\lbrace 0\rbrace$ for $k<6$.
We shall prove
\bigskip
\noindent{\bf Theorem.} {\it Let $\sigma\in {\bf R}$ be fixed. Then the following assertions hold:
\medskip
i) If $k\geq 4$, then there is a Hecke eigenform $f\in S_{k+{1\over 2}}^{(+)}(4)$ with $L^*(f,\sigma)\neq 0$.

ii) If $k\geq 6$ and $\sigma\neq {k\over 2}+{1\over 4}$, then there is a Hecke eigenform $f\in S_{k+{1\over 2}}^{(-)}(4)$ with $L^*(f,\sigma)\neq 0$.}
\bigskip
\noindent{\it Remarks. i)} In the proof of the theorem we explicitly construct a cusp form $f$ in the space $S_{k+{1\over 2}}^{(+)}(4)$ (resp. in the space $S_{k+{1\over 2}}^{(-)}(4)$) such that $L^*(f,\sigma)\neq 0$, for any real $\sigma$. 

\noindent {\it ii)}  Let $k=4$. It follows from sect. 2 that $S_{9/2}^{(+)}(4)={\bf C}\theta\Delta_4$, hence $f_1:=\theta\Delta_4$ is a Hecke eigenform and $L^*(f_1,\sigma)\neq 0$ for all $\sigma\in {\bf R}$. It will in fact follow from the proof given in the next sect. that $L^*(f_1,\sigma)> 0$ for all $\sigma\in {\bf R}$. Similarly, in the case $k=6$ one checks that $S_{13/2}^{(-)}(4)={\bf C}\theta D_2\Delta_4$, so  $f_2:=\theta D_2\Delta_4$ is a Hecke eigenform and $L^*(f_2,\sigma)\neq 0$ for all $\sigma\in {\bf R}, \sigma\neq {{13}\over 4}$. It will follow from our proof that $L^*(f_2,\sigma)$ is positive for $\sigma > {{13}\over 4}$ and negative for  $\sigma < {{13}\over 4}$.\bigskip
\noindent{\bf 4. Proofs}
\bigskip
We start with the proof of i) whose first part is standard. If $f\in S_{k+{1\over 2}}(4)$, then in the usual way by Mellin transform one has
$$L^*(f,s)=\int_0^\infty f(it)(2t)^s{{dt}\over t} \quad (\sigma \gg 1).$$
Splitting up the integral into the sum of the integral from 0 to ${1\over 2}$ and the integral from ${1\over 2}$ to $\infty$, and substituting $t\mapsto {1\over {4t}}$ in the former one, we obtain
$$L^*(f,s)=\int_{1\over 2}^\infty f({i\over {4t}})({1\over {2t}})^s{{dt}\over t} \,+ \, \int_{1\over 2}^\infty f(it)(2t)^s{{dt}\over t}.$$
If $f|W_4=f$, then $f({i\over {4t}})=(2t)^{k+1/2}f(it)$, hence we get the identity
$$L^*(f,s)=\int_{1\over 2}^\infty f(it)\Bigl((2t)^s+(2t)^{k+1/2-s}\Bigr){{dt}\over t} \leqno (8)$$
which is valid for all $s\in {\bf C}$, since $f$ is a cusp form.
\medskip
We now put $f:=\Delta_4\theta^{2k-7}$. By (4) and (6) we see that $f\in S_{k+{1\over 2}}^{(+ )}(4)$.
\medskip
Recall that the valence formula asserts that the sum of the orders (counted with multiplicities) of a non-zero modular form $g$ of weight $\ell$ on $\Gamma_0(4)$ on the compact modular curve $X_0(4)$ of level 4 is equal to ${{\ell}\over {12}}[\Gamma_1:\Gamma_0(4)]$. We note  that $[\Gamma_1:\Gamma_0(4)]=6$ and apply the valence formula with $g=\Delta_4\in M_4(4)$. Since $\Delta_4$ vanishes at the cusps 0 and $i\infty$, it follows that $\Delta_4$ has no zeros on $\cal H$, in particular $\Delta_4(it)\neq 0$ for all $t>0$. Since by definition $\Delta_4(it)$ is real for $t>0$, by continuity we must have $\Delta_4(it)>0$ for all $t$ or $\Delta_4(it)<0$ for all $t$. We claim that the first alternative holds. Indeed, evaluating (5) at $z={i\over 2}$ we find that
$$\theta^4({i\over 2})=32F_2({i\over 2}),$$ hence
$$\Delta_4({i\over 2})=F_2({i\over 2})\Bigl(\theta^4({i\over 2})-16 F_2({i\over 2})\Bigr)>0.$$
(Alternatively, we could also have used the convergent $q$-product expansion
$$\Delta_4(z)=q\prod_{n\geq 1,n\equiv {0, \pm 1}\pmod 4}(1-q^n)^8,$$
cf. [2, p. 25].)

By definition, $\theta(z)$ is real and positive on the positive imaginary axis. 
Thus $f(it)>0$ for $t>0$ and so in particular the integrand in (8) is positive for all $t\geq {1\over 2}$ if $s=\sigma$ is real, hence $L^*(f,\sigma)$ is positive and so non-zero for all $\sigma$. Writing $f$ in terms of a Hecke basis of $S_{k+{1\over 2}}^{(+ )}(4)$ we deduce the assertion of i).
\medskip
We now prove ii). If $f|W_4=-f$, then similarly as in i) we find that
$$L^*(f,s)=\int_{1\over 2}^\infty f(it)\Bigl((2t)^s-(2t)^{k+1/2-s}\Bigr){{dt}\over t}.$$
We let $f:=\Delta_4D_2\theta^{2k-11}$. By (4), (6) and (7) then $f\in S_{k+{1\over 2}}^{(- )}(4)$.
\medskip
We want to show that $L^*(f,\sigma)\neq 0$ for $\sigma$ real, $\sigma\neq {k\over 2}+{1\over 4}$. By the functional equation (3) we may assume that $\sigma >{k\over 2}+{1\over 4}$. In this range clearly
$$(2t)^\sigma-(2t)^{k+1/2-\sigma}\geq 0$$
whenever $t\geq {1\over 2}$, and the inequality is strict for $t>{1\over 2}$.
\medskip
As mentioned above, $\theta(it)>0$ and we proved that $\Delta_4(it)>0$, for $t>0$. Thus to prove our claim, it will suffice to show that $D_2(it)>0$ for $t>{1\over 2}$.
\medskip
We apply  the valence formula with $g=D_2\in M_2(4)$. Observe  that $D_2(z)$ by (7) vanishes at $z={i\over 2}$. We conclude that any zero of $D_2$ in ${\cal H}$ must be equivalent under  $\Gamma_0(4)$ to ${i\over 2}$. Thus to see that $D_2(it)\neq 0$ for $t>{1\over 2}$ we have to show that for $t>{1\over 2}$ the point $it$ cannot be equivalent to ${i\over 2}$ under $\Gamma_0(4)$. \medskip Assume on the contrary that there exists $M\in\Gamma_0(4)$ with $M\circ {i\over 2}=it$ and $t>{1\over 2}$. Let ${\cal F}$ be the standard fundamental domain for $\Gamma_1$ consisting of $z\in {\cal H}$ with $|z|\geq 1$ and $|x|\leq {1\over 2}$. Let $S=\pmatrix {0&{-1}\cr 1&0\cr}$. Then $S^2=-E$. Suppose first that $t\leq 1$. Then $S\circ it\in {\cal F}$. Also $S\circ {i\over 2}=2i\in int {\cal F}$ and $SMS\circ 2i=S\circ it$. It follows that $SMS=\pm E$  and so $2i={i\over t}$, i.e. $t={1\over 2}$, a contradiction.\medskip
Now suppose that $t>1$. Then $M\circ {i\over 2}\in {\cal F}$. Since $M\circ {i\over 2}=MS\circ (S\circ {i\over 2})=MS\circ 2i$ and $2i\in int{\cal F}$, it follows that $MS=\pm E$, i.e. $M=\pm S$, a contradiction since $S$ is not in $\Gamma_0(4)$.
\medskip
Since $D_2$ has the value 1 at infinity, by continuity we finally find that indeed $D_2(it)>0$ for all $t>{1\over 2}$.
\medskip
It follows that $L^*(f,\sigma)\neq 0$ as claimed. Writing $f$ in terms of a Hecke basis as before, we derive the assertion of ii).

\bigskip \noindent {\it Acknowledgements.}   The first author was partially supported by NRF 2017R1A2B2001807.
\bigskip\bigskip

\bigskip
\noindent{\bf References}
\bigskip
\noindent [1] H. Cohen, Sums involving the values at negative integers of $L$-functions of quadratic 
characters, Math. Ann. 217  (1975), no. 3, 271-285
\medskip
\noindent [2]  H. Cohen, Formes modulaires \` a une et deux variables, Th\`ese, Univ. Bordeaux I, 1976
\medskip
\noindent [3] S. Das and W. Kohnen, Non-vanishing of Koecher-Maass series attached to cusp forms,
Adv. Math. 281 (2015), 624-669
\medskip
\noindent [4] W. Kohnen, Nonvanishing of Hecke $L$-functions associated to cusp forms inside the
critical strip, J. of Number Theory 67 (1977), vol. 2, 182-189
\medskip
\noindent [5] W. Kohnen, Non-vanishing of Koecher-Maass series attached to Siegel cusp forms on
the real line, Abh. Math. Sem. Univ. Hamburg  87 (2017), no. 1, 39-41
\medskip
\noindent [6] M. Knopp, W. Kohnen and W. Pribitkin, On the signs of Fourier coefficients of cusp forms,
Ramanujan J. 7 (2003), no. 1-3, 269-277
\medskip
\noindent [7]  A. Raghuram, Non-vanishing of $L$-functions of cusp forms inside the critical strip, 
Number theory, 97-105, Ramanujan Math. Soc. Lect. Notes Ser. 1, Ramanujan Math. Soc., Mysore (2005)
\medskip
\noindent [8] B. Ramakrishnan and Karam Deo Shankhadhar, Non-vanishing of $L$-functions associated
to cusp forms of half-integral weight, Automorphic Forms, Springer Proc. Math. Stat. 115 (2014), 223-231
\medskip
\noindent [9] G. Shimura, On modular forms of half-integral weight, Ann. of Math. 97 (1973), 440-481
\medskip
\noindent [10] H. Yoshida, On calculations of zeros of various $L$-functions, J. Math. Kyoto Univ. 35, no. 4 (1995), 663-696

\bigskip\bigskip
\noindent {\it Department of Mathematics, Pohang Institute of Science and Technology, POSTECH, Pohang 790-784, Korea}
\medskip
\noindent {\it E-mail: yjc@postech.ac.kr}\bigskip
\noindent {\it Mathematisches Institut der Universit\"at, INF 205, D-69120 Heidelberg, Germany}
\medskip
\noindent{\it E-mail: winfried@mathi.uni-heidelberg.de}

\bye